\documentclass[11pt,a4paper]{article}
\usepackage{t1enc}
\usepackage[latin1]{inputenc}
\usepackage[english]{babel}
\begin{document}

\title{Remark on Periodic Solutions of Non Linear Oscillators}
\author{A. Raouf  Chouikha \bigskip \\ University of Paris-Nord \\ Institut Galilee \\ LAGA,CNRS UMR 7539\\ 
 Villetaneuse F-93430 \\ \small{Chouikha@math.univ-paris13.fr}}
\date{}

\maketitle

\begin{abstract}
We contribute to the method of trigonometric series for solving differential equations of certain non linear oscillators.\footnote{ 2000 {\it Mathematical Subjects Classification,} 34A20. }

\end{abstract}

 \bigskip

\section{Introduction}
The non linear nonharmonic motion of an oscillator may be given by the following differential equation 
\begin{equation}
u'' + \omega^2 u = - \beta u^2
\end{equation}

$\beta $ \ and\ $\omega $\ being constants, with initial conditions 
\begin{equation}
u(0) = a_0, \qquad u'(0) = 0
\end{equation}

To solve this problem A. Shidfar and A. Sadeghi [1], have given two series solutions. They describe a general approach in which the differential equation, rather than the solutions series , is majorized.\\ Notice that if we write $a_0 = - {\omega^2 \over \beta} $ then \ 
\begin{equation}
u(t) \equiv - {\omega^2 \over \beta}
\end{equation}
 is a trivial solution of (1) and (2). \\ They gave a series solutions of (1) and (2), which includes (3) as a special case.\\

By writing $$ u = v - {\omega^2 \over {2\beta}}$$ the problem becomes 
\begin{equation}
v'' + \beta v^2=  {\omega^4 \over {4\beta}}
\end{equation}
under the initial conditions
\begin{eqnarray} \cases{ 
v(0) = a_0 + {\omega^2 \over {2\beta}} & \cr
v'(0) = 0 & \cr}
\end{eqnarray}
The method of [1] consists to solve equations (4) and (5) in the form 
\begin{equation}
v(t) = c_0 + c_1 \sin\omega t + c_2 \sin^2 \omega t + c_3 \sin^3 \omega t  + .....
\end{equation}
where  \ $ c_i , i = 0, 1, 2, .... $\ are coefficients to be determined by the substitution of (6) in (4).\\
In fact, \  $\omega = {\pi \over T}$ \ where $T$ is the  period of the solution, which can be expressed in terms of the Weierstrass function \ $\wp(z,2T , 2T ')$.\\ 
So, we find that $$ 2 \omega^2 c_2 + \beta c_0^2 = {\omega^2 \over {4\beta}}.$$
For $ n  \geq 1$,  the recursion formula for these  coefficients are
\begin{equation}
(n+1)(n+2) c_{n+2} = n^2 c_n - {\beta \over \omega^2} \sum_{r=0}^{n} c_r c_{n-r}.
\end{equation}

Equations (5) and (6) imply that \ $c_1 = 0$.\ Relations (7) yields  $$c_3 = 0, \quad c_5 = 0, .....$$

The even order coefficients simply are 
$$
c_0 = a_0 + {\omega^2 \over {2\beta}}, $$ 
$$c_2 = - {a_0 \over{2\omega^2}}(\omega^2 + a_0 \beta), $$
$$ c_4 =  - {\beta \over{6\omega^2}}a_0 (\omega^2 + a_0 \beta)({3\over 4} - {a_0\beta \over{2\omega^2}}),  $$
$$c_6 =  - {\beta \over{180\omega^2}}a_0 (\omega^2 + a_0 \beta)({3\over 4} - {a_0\beta \over{2\omega^2}})(15 - {2 a_0\beta \over{\omega^2}}) - {\beta a_0^2 \over{120\omega^6}}(\omega^2 + a_0 \beta)^2,$$ \quad {\it etc.}\\
The coefficient \ $c_0$\ follows from the condition (5). The solution for the equations (4) and (5) can now be written as

\begin{equation}
u(t) = a_0 - {a_0 \over{2\omega^2}}(\omega^2 + a_0 \beta) \sin^2 \omega t + ....
\end {equation}
Relations of the coefficients and further induction show that \quad $c_{2i}, \ i=1,2,...$\quad all vanish for \quad $a_0 =  {\omega^2 \over {\beta}} .$ \quad  So, the trivial solution (3) is included in (6) as a special case. 

\section{Convergence of the solutions}

We now show the convergence of these series. In [1] one proved the following 
\bigskip

{\bf Lemma 1}  \qquad {\it The serie (6) solution of Equation (4)-(5) is absolutely convergent for all} $t$.\\

{\bf Proof} \quad We firstly note that if \ $c_0 > 0,\ c_2 > 0$\ and\ $\beta < 0,$ \ then all coefficients \ $c_n$\ in the serie expansion (6) are positive. Indeed, we may write $$\sum_{n\geq 0} ( (n+1)(n+2) c_{n+2} = \sum_{n\geq 0} n^2 c_n - {\beta \over{\omega^2}}(\sum_{n\geq 0}  c_n )^2 + {\beta \over{\omega^2}} c_0^2,$$
or 
\begin{equation}
- \beta (\sum_{n\geq 0}  c_n )^2 + {\omega^2} \sum_{n\geq 0} n c_n = - \beta c_0^2 - 2 {\omega^2} c_2.
\end{equation}

Since the right hand side of (9) is finite and \ $c_i$\ are positive, the series \quad $\sum_{n\geq 0} c_n$\quad converges. Following  [1], if we put $$c'_0 = \mid c_0\mid, \quad c'_1 = 0, \quad c'_2 = \mid c_2\mid $$
and for \ $n \geq 2$
 $$c'_{n+2} = {n^2\over{(n+1)(n+2)}} c'_n  + {\mid \beta \mid \over {\omega^2(n+1)(n+2)}} \sum_{r=0}^{n} c'_r c'_{n-r},$$
then the series \ $\sum_{n\geq 0} c'_n$ \ converges. Since \ $\mid c_n\mid \geq c'_n$, \ it follows that the solution series (6) is absolutely convergent, and hence the series expansion solution of (1)-(2) converges for all $t. $ \\

We notice that  we may deduce Lemma 1 a previous result concerning  Equation (4).\\ We have shown that the coefficients verify a more general properties. Indeed, we have [4]

\bigskip

{\bf Lemma 2}\qquad  {\it For any positive number \ $\epsilon $ \ small enough (but $\epsilon  \neq 0$), there exists a positive constant $k$  \ verifying $$ k < {\beta \over \omega^2}{3\over 4}\epsilon 4^{\epsilon -{1\over 2}}$$
such that the coefficients $c_n$ of the series expansion (6) solution of the differential equation  (4)-(5) satisfy the inequality
\begin{equation}
\mid c_n \mid < {k \over{n^{{3\over 2}-\epsilon }}}.
\end{equation} }
\\ 

{\bf Proof}\quad We first notice that Lemma 2 gives an optimal result, because our method do not run for \ $\epsilon  = 0$.\\ 
The coefficients \ $c_n$\ of the power series  solution, satisfy the recursion formula (7). We shall prove there exist two positive constants \ $k > 0,$ and $\alpha > 1$, such that  the following inequality holds 
$$\mid c_n \mid < {k\over{n^\alpha }}$$ for any integer $n \geq 1$.   
Suppose for any \ $n \leq p$\ , we get $\mid c_n \mid < {k\over{n^\alpha }}$. 	In particular, it implies that 
$$\sum_{0<r<p} c_r c_{p-r} < \sum_{0<r<p} {k^2\over{r^\alpha (r - p)^\alpha }} \leq {k^2\over{(p-1)^{\alpha -1}}}.$$ Equality (7) gives $$ c_{p+2} = {p^2 - 2{\beta \over \omega^2}c_0\over{(p+1)(p+2)}} c_p - {\beta \over \omega^2 (p+1)(p+2)} \sum_{r=1}^{p-1} c_r c_{p-r}.$$ 
Thus, if we prove the following inequality 
\begin{equation}
{{p^2-2{\beta \over \omega^2}c_0} \over{(p+1)(p+2)}}{k\over{p^\alpha }}  + {\beta \over \omega^2 (p+1)(p+2)}{k^2\over{(p-1)^{\alpha -1}}} \leq {k\over{(p+2)^\alpha }}
\end{equation}
 so   
\begin{equation}
\mid c_{p+2} \mid < {k\over{(p+2)^\alpha }}
\end{equation}
Notice that (11) implies that \ $$ {{p^2-2{\beta \over \omega^2}c_0} \over{(p+1)(p+2)}}{k\over{p^\alpha }}   \leq {k\over{(p+2)^\alpha }}.$$
Thus,  a necessary condition to (11) holds is : \ $\alpha \leq {3\over 2}.$\  
Inequality (11) is equivalent to 
\begin{equation}
k < {\beta \over \omega^2} p f(p) g(p)
\end{equation}
where $$f(p) = {p+1\over{p}} ({p-1\over{p+2}})^{\alpha -1}$$
 $$g(p) = 1 - {(p^2 - {\beta \over \omega^2}c_0)(p+2)^{\alpha -1}\over{(p+1)p^\alpha }}$$
By using MAPLE, we are able to prove that $f(p)$ is an increasing positive function in $p$ . Moreover, for any $p \geq 1$, $f(p)$ is minorated  $$f(p) \geq ({3\over 2}) 4^{1-\alpha }.$$ The function $g(p)$ is such that   $$p g(p) = p - {(p^2 - {\beta \over \omega^2}c_0)({p+2\over p})^{\alpha -1}\over{(p+1)}}$$ is a  strictly decreasing and bounded function .\\ More exactly, we may calculate the lower bound
$$ g(p) >{ (3 - 2\alpha )\over p}.$$
Thus, if \ $(3 - 2\alpha ) = \epsilon > 0,$\ it suffices to choice $$k \leq  ({3\over 2}) 4^{1-\alpha }(3 - 2\alpha )$$ to inequality (13) holds.  \\

{\bf Remark for the case $\epsilon = 0$ :}\\
Notice that the choice of $k$ depends on $\alpha $ value.\\
For $\alpha  = {3\over 2},$ we then prove by MAPLE that  the function $$p g(p) = p - {(p^2 - { 3\over 2\omega^2}c_0)({p+2\over p})^{1\over 2}\over{(p+1)}}$$ is positive and strictly decreasing  to $0$. While $p^2 g(p)$ is a bounded function. \\ Moreover, it appears that  $p f(p) g(p)$ is a decreasing function which tends to $0$ when $p$ tends to infinity. Thus, our method falls since it do not permit to determine a non negative constant $k$.\\

\bigskip

Following [1], it is interesting to write the power series solution for the system (4)-(5), 
\begin{equation}
v(x) = \sum_{n=0} ^\infty  b_n x^n.
\end{equation}
  We find again that \ $$b_{2p+1} = 0 \qquad p = 0,1,2,...,$$ while $$b_0 = a_0 +  {\omega^2 \over {2\beta}},$$
$$2 b_2 + \beta  b_0^2 = {\omega^4 \over {4\beta}},$$
$$(n+2)(n+1)b_{n+2} = -\beta  \sum _{r=0}^{r=n} b_r b_{n-r},$$ where  $n$ is even and non zero.\\
The coefficients\ $b_{2p}, p = 0,1,2,...,$\ again vanish for $a_0 = - {\omega^2 \over {\beta}}.$ \ We may verify that the solutions  (6) and (14) are identical.\\ The latter method permits  to compare approximate solutions of the anharmonic motion of the oscillator.\\

\section{Another differential equation}

We now examine the following differential equation

\begin{equation}
u'' + \omega^2 u = - \beta u^3
\end{equation}

$\beta $ \ and\ $\omega $\ being constants, with initial conditions 
\begin{equation}
u(0) = a_0, \qquad u'(0) = 0.
\end{equation}

We put \ $v = {u\over a_ 0}$\ and \ $t = \omega x.$\
 We then obtain from (15) and (16) 
\begin{equation}
{d^2 v\over{dt^2}} + v + \beta v^3 = 0, \quad v(0) = 1, \quad v'(0) = 0 
\end{equation}
where $ \beta = {\beta a_0^2\over \omega }.$\\
A. Shidfar and A. Sadeghi [2] solved (17) by series method in Sinus power
\begin{equation}
v(t) = c_0 + c_1 \sin\omega t + c_2 \sin^2 \omega t + c_3 \sin^3 \omega t  + .....
\end{equation}
Here, \  $\omega = {\pi \over T}$ \ where $T$ is the  period of the solution, which can be expressed in terms of the Jacobi function \ $sn(z,2T , 2T')$.\\
So, $$c_0 = a_0.$$
For $n \geq 1$, we get the recursion formula 
\begin{equation}
(n+1)(n+2) c_{n+2} = n^2 c_n - {\beta \over \omega^2} \sum_{r=0}^{n} \sum_{m=0}^{n-r}c_m c_r c_{n-m-r}.
\end{equation}

 Under some conditions, they proved  estimates of the coefficients 
$$\mid c_n \mid \leq R^n$$ where \ ${1\over R}$\ is a radius of convergence. \\
In fact, we may prove an analog of Lemma 2 for this equation. Indeed, we have 
$$ \sum_{r=0}^{n} \sum_{m=0}^{n-r}c_m c_r c_{n-m-r} = 2 c_0^2 c_n + 2 c_0 c_1c_{n-1} +c_0 \sum_{m=0}^n c_n c_{n-m} + \sum_{r=2}^{n-2} \sum_{m=0}^{n-r}c_m c_r c_{n-m-r}.$$
Then, $$\sum_{r=2}^{n-2} \sum_{m=0}^{n-r}c_m c_r c_{n-m-r} = \sum_{r=2}^{n-2} c_r [2 c_0 c_{n-r} +  \sum_{m=1}^{n-r-1}c_m  c_{n-m-r}] $$ 
$$< \sum_{r=2}^{n-2} \mid c_r \mid  [2\mid c_0 c_{n-r} \mid + {k^2\over{(n-r-1)^{\alpha-1}}}]$$
  $$< 2 c_0 k^2 \sum_{r=2}^{n-2} {1\over{r^\alpha (n-r)^\alpha}}  + k^3 \sum_{r=2}^{n-2} {1\over{r^\alpha (n-r-1)^{\alpha-1}}} < {2 c_0 k^2\over{(n-1)^{\alpha-1}}} + {k^3\over{(n-2)^{\alpha -1}}}$$
 Finally,  $$ \mid c_{p+2} \mid < {k\over{(p+2)^\alpha }}$$
 as soon as the non negative constant $k$ satisfies the inequality
$$ {k\over{n^{\alpha -2}}} + {2 c_0 k^2\over{n^\alpha}} + {2 c_0 c_1k\over{(n-1)^{\alpha }}} + {3 c_0 k^2\over{(n-2)^{\alpha-1}}} + {k^3\over{(n-2)^{\alpha -1}}} < {(n+1) k \over{(n+2)^{\alpha -1}}},$$ So, 
$$ {1\over{n^{\alpha -2}}} + {k^2 +3 c_0 k\over{(n-2)^{\alpha-1}}} + {2 c_0 + 2 c_0 c_1\over{(n-1)^{\alpha }}}  < {(n+1) \over{(n+2)^{\alpha -1}}}.$$
 By using MAPLE, we verify it is possible to find a such constant. Moreover, we find again the necessary condition : \
 $\alpha  \leq {3\over 2},$\ since we get \ ${1\over{n^{\alpha -2}}} < {(n+1) \over{(n+2)^{\alpha -1}}}.$\\

{\bf General remarks :} \qquad It is wellknown from the theory of elliptic functions that solutions of equations  (4) and (17) are related. This allows to express the series expansion of the solution of (17) from a series expansion of a solution of (4) and conversely. Indeed, one has $$\wp(z)  = C - {\delta^2 \over{sn^2 (\delta z)}},$$ where\ $\wp(z)$\ is the elliptic Weierstrass function  and\ $sn(u)$\ is the Jacobi fonction , $\delta $ is a constant , only dependent on the initial parameters. \\ Notice that series expansion of the\ $sn(u)$ \ in sinus power was given in a previous paper  (see Proposition (2.1) in [3]).\\

\bigskip
\begin{center}
{\bf REFERENCES } 
\end{center}

\bigskip
[1] \ A. Shidfar and A. Sadeghi \quad {\it Some Series Solutions of the Anharmonic Motion Equation}, \quad J. Math. Anal. Appl.  120, p. 488-493 (1986).

[2] \ A. Shidfar and A. Sadeghi \quad {\it The Periodic Solutions of Certain Non-linear Oscillators},\quad Appl. Math. Lett., vol 3, n 4, p. 21-24, (1990).  

[3] \ R. Chouikha \quad {\it Sur des developpements de fonctions elliptiques}, \quad Publ. Math. Fac des Sci de Besancon, Fasc. Th des Nombres, p.1-9,  (1989).

[4] \ R. Chouikha \quad {Fonctions elliptiques et bifurcations d'equations differentielles}, \quad Canad.  Math. Bull.,  vol. 40 (3), p. 276-284, (1997)

 \end{document}